\overfullrule=0pt
\centerline {\bf Singular points of non-monotone potential operators}\par
\bigskip
\bigskip
\centerline {BIAGIO RICCERI}\par
\bigskip
\bigskip
\centerline {\it Dedicated to the memory of Francesco Saverio De Blasi}
\bigskip
\bigskip
{\bf Abstract:} In this paper, we establish some results about the singular points of certain non-monotone potential operators.
Here is a sample: If $X$ is an infinite-dimensional reflexive real Banach space and if $T:X\to X^*$
is a non-monotone, closed, continuous potential operator such that the functional $x\mapsto \int_0^1T(sx)(x)ds$
is sequentially weakly lower semicontinuous and $\lim_{\|x\|\to +\infty}(\int_0^1T(sx)(x)ds+\varphi(x))=+\infty$
for all $\varphi\in X^*$, then the set of all singular points of $T$ is not $\sigma$-compact.\par
\bigskip
\bigskip
\bigskip
\bigskip
Here and in what follows, $(X,\|\cdot\|)$ is a reflexive real
Banach space, with topological dual $X^*$, and $T:X\to X^*$ is a continuous potential operator. 
This means that  the functional 
$$x\to J_T(x):=\int_0^1T(sx)(x)ds$$
is of class $C^1$ in $X$ and its G\^ateaux derivative is equal to $T$.
\smallskip
Let us recall a few classical definitions.\par 
\smallskip
Def. 1. - $T$ is said to be monotone if
$$(T(x)-T(y))(x-y)\geq 0$$
for all $x, y\in X$.\par
\smallskip
  This is equivalent to the fact that the functional $J_T$ is convex.\par
\smallskip
Def. 2. - $T$ is said to be closed if for each closed set $C\subseteq X$, the set $T(C)$ is closed in $X^*$.\par
\smallskip
Def. 3. - $T$ is said to be compact if for each bounded set $B\subset X$, the set $\overline {T(B)}$ is compact in $X^*$.\par
\smallskip
Def. 4. - $T$ is said to be proper if for each compact set $K\subset X^*$, the set $T^{-1}(K)$ is compact in $X$.\par
\smallskip
Def. 5. - $T$ is said to be a local homeomorphism at
a point $x_0\in X$ if there are a neighbourhood $U$ of $x_0$ and a neighbourhood
$V$ of $T(x_0)$ such that the restriction of $T$ to $U$ is a homeomorphism between $U$
and $V$. If $T$ is not a local homeomorphism at $x_0$, we say that $x_0$ is
a singular
point of $T$.\par
\smallskip
We denote by $S_T$ the set of all singular points of $T$. Clearly, the set $S_T$ is closed.\par
\smallskip
Assume that the restriction of $T$ to some open set $A\subseteq X$ is of
class $C^1$.\par
\smallskip
 We  then denote by
 $\tilde S_{T_{|A}}$ the set of all $x_0\in A$ such that
the operator
$T'(x_0)$ is not invertible. Since the set of all invertible operators belonging to ${\cal L}(X,X^*)$ is
open in ${\cal L}(X,X^*)$, by the continuity of $T'$,
the set $\tilde S_{T_{|A}}$ is closed in $A$.\par
\smallskip
Def. 6. -  $T$ is said to be a Fredholm operator of index zero in $A$
if, for each $x\in A$, the codimension of $T'(x)(X)$ and  the dimension of
$(T'(x))^{-1}(0)$ are finite and equal.
\smallskip
Def. 7. - A set in a topological space is said to be $\sigma$-compact if it is the union of an at most countable
family of compact sets.\par
\smallskip
Def. 8. - A functional $I:X\to {\bf R}$ is said to be coercive if
$$\lim_{\|x\|\to +\infty}I(x)=+\infty\ .$$
The aim of this note is to establish the  following results:\par
\medskip
THEOREM 1. - {\it If $X$ is infinite-dimensional, $T$ is closed and non-monotone,  $J_T$ is
 sequentially weakly lower
semicontinuous and $J_T+\varphi$ is coercive
for all $\varphi\in X^{*}$, 
then 
both $S_T$ and $T(S_T)$ are not $\sigma$-compact.}\par
\medskip
THEOREM 2. - {\it In addition to the assumptions of Theorem 1, suppose that
there exists a
closed, $\sigma$-compact set $B\subset X$ such that the restriction
of $T$ to $X\setminus B$ is of class $C^1$.\par
Then, both $\tilde S_{T_{|(X\setminus B)}}$ and
$T(\tilde S_{T_{|(X\setminus B)}})$ are not $\sigma$-compact.}\par
\medskip
THEOREM 3. - {\it Assume that $(X,\langle\cdot,\cdot\rangle)$ is a Hilbert space, with $\hbox {\rm dim}(X)\geq 3$, and that $T$ is compact and
of class $C^1$ with
$$\liminf_{\|x\|\to +\infty}{{J_T(x)}\over {\|x\|^2}}\geq 0 \eqno{(1)}$$ 
and, for some $\lambda_0\geq 0$,
$$\lim_{\|x\|\to +\infty}\|x+\lambda T(x)\|=+\infty \eqno{(2)}$$
for all $\lambda>\lambda_0$ .\par
Set
$$\Gamma=\{(x,y)\in X\times X : \langle T'(x)(y), y\rangle<0\}$$
and, for each $\mu\in {\bf R}$, 
$$A_{\mu}=\{x\in X :  T'(x)(y)=\mu y\hskip 5pt for\hskip 5pt some\hskip 5pt y\in X\setminus \{0\}\}\ .$$
When $\Gamma\neq\emptyset$, set also
$$\tilde\mu=\max\left \{ -{{1}\over {\lambda_0}}, \inf_{(x,y)\in \Gamma}{{\langle T'(x)(y), y\rangle}\over {\|y\|^2}}\right
\}\ .$$
Then, the following assertions are equivalent:\par
\noindent
$(i)$\hskip 5pt the operator $T$ is not monotone\ ;\par
\noindent
$(ii)$\hskip 5pt there exists $\mu<0$ such that $A_{\mu}\neq\emptyset$\ ;\par
\noindent
$(iii)$\hskip 5pt $\Gamma\neq\emptyset$ and, for each $\mu\in ]\tilde\mu,0[$, the set
$A_{\mu}$ contains an accumulation point\ .}\par
\medskip
REMARK 1. - Of course, Theorem 2 is meaningful only when $X$ and $X^*$ are linearly isomorphic. Indeed, if not, 
the fact that $\tilde S_{T_{|(X\setminus B)}}$  is not $\sigma$-compact follows directly from the equality
 $\tilde S_{T_{|(X\setminus B)}}=X\setminus B$ \ .\par
\medskip
The previous theorems extend and improve the results of [3] in a remarkable way. The reason for this
resides in the tools used to prove them. Precisely, in [3], the main tools were Theorems A and B below
 jointly with the minimax theorem proved in [2]. This latter theorem contains a severe restriction: one of
the two variables on which the underlying function depends must run over a real interval. In the current paper,
we still continue to use  Theorems A and B in an essential way but, this time, jointly with a consequence of another very recent
minimax theorem ([4], Theorem 3.2) which is not affected by the above recalled restriction.\par
\medskip
So, let us recall Theorems A and B.\par
\medskip
THEOREM A ([6], Theorem 2.1). - {\it If $X$ is infinite-dimensional, if  $T$
is closed and if $S_T$ is $\sigma$-compact, then the
restriction of $T$ to $X\setminus S_T$ is a homeomorphism
between $X\setminus S_T$ and $X^*\setminus T(S_T)$.}\par
\medskip
THEOREM B ([1], Theorem 5). - {\it If \hbox {\rm dim}$(X)\geq 3$,
if $T$ is a $C^1$
proper Fredholm operator of index zero and if $\tilde S_T$ is discrete,
then $T$ is a homeomorphism between $X$ and $X^*$.}\par
\medskip
As we said above, besides Theorems A and B, the other major tool that we will use is a consequence of the following minimax
theorem (here stated in a particular version which is enough for our purposes):\par
\medskip
THEOREM C ([4], Theorem 3.2). - {\it Let $Y$ be a convex set in a real vector space $E$ and
let $f:X\times E\to {\bf R}$ be sequentially weakly lower semicontinuous and coercive in $X$,
and linear in $E$. Assume also that
$$\sup_Y\inf_Xf<\inf_X\sup_Yf\ .$$
Then there exists $\tilde y\in Y$ such that the functional $f(\cdot,\tilde y)$ has at least two global minima.}\par
\medskip
\par
Let us introduce the following notations. We denote by ${\bf R}^X$ the space of
all functionals $\varphi:X\to {\bf R}$.  For each $I\in {\bf R}^X$ and for
each  of non-empty subset $A$ of $X$, we denote by
$E_{I,A}$ the set of all $\varphi\in {\bf R}^X$ such that $I+\varphi$ is sequentially
weakly lower semicontinuous and coercive, and
$$\inf_A\varphi\leq 0\ .$$
Here is the above mentioned consequence of Theorem C:\par
\medskip
THEOREM 4. - {\it Let  $I:X\to {\bf R}$ be a functional and $A, B$ two non-empty subsets of $X$ such that
$$\sup_AI<\inf_BI\ .\eqno{(3)}$$
Then, for every convex set $Y\subseteq E_{I,A}$ such that
$$\inf_{x\in B}\sup_{\varphi\in Y}\varphi(x)\geq 0\hskip 5pt and\hskip 5pt  
\inf_{x\in X\setminus B}\sup_{\varphi\in Y}\varphi(x)=+\infty\ ,\eqno{(4)}$$
 there exists $\tilde\varphi\in Y$ such that the
functional $I+\tilde\varphi$ has at least two global minima.}\par
\smallskip
PROOF. 
Consider the function $f:X\times {\bf R}^X\to {\bf R}$ defined by
$$f(x,\varphi)=I(x)+\varphi(x)$$
for all $x\in X$, $\varphi\in {\bf R}^X$. Fix $\varphi\in Y$. In view of $(3)$, we also can fix $\epsilon\in ]0,\inf_BI-\sup_AI[$.
Since $\inf_A\varphi\leq 0$, there is $\bar x\in A$ such that $\varphi(\bar x)<\epsilon$. Hence,
we have
$$\inf_{x\in X}(I(x)+\varphi(x))\leq I(\bar x)+\varphi(\bar x)<\sup_AI+\epsilon\ ,$$
from which it follows that
$$\sup_{\varphi\in Y}\inf_{x\in X}(I(x)+\varphi(x))\leq \sup_AI+\epsilon<\inf_BI\ .\eqno{(5)}$$
On the other hand, in view of $(4)$, one has
$$\inf_BI\leq \inf_{x\in B}(I(x)+\sup_{\varphi\in Y}\varphi(x))=\inf_{x\in B}\sup_{\varphi\in Y}(I(x)+\varphi(x))=
\inf_{x\in X}\sup_{\varphi\in Y}(I(x)+\varphi(x))\ .\eqno{(6)}$$
Finally, from $(5)$ and $(6)$, it follows that
$$\sup_{\varphi\in Y}\inf_{x\in X}f(x,\varphi)<\inf_{x\in X}\sup_{\varphi\in Y}f(x,\varphi)\ .$$
Therefore the function $f$ satisfies the assumptions of Theorem C and the conclusion follows.\hfill $\bigtriangleup$
\medskip
More precisely, we will use the following corollary of Theorem 4:\par
\medskip
COROLLARY 1. - {\it Let $I:X\to {\bf R}$ be a sequentially weakly lower semicontinuous, non-convex functional such that
$I+\varphi$ is coercive for all $\varphi\in X^*$.\par
Then, for every convex set $Y\subseteq X^*$ dense in $X^*$, there exists $\tilde\varphi\in Y$ such that the functional
$I+\tilde\varphi$ has at least two global minima.}\par
\smallskip
PROOF. Since $I$ is not convex,  there
exist $x_1, x_2\in X$ and $\lambda\in ]0,1[$ such that
$$\lambda I(x_1)+(1-\lambda) I(x_2)<I(x_3)$$
where
$$x_3=\lambda x_1+(1-\lambda)x_2\ .$$
Fix $\psi\in X^*$ so that
$$\psi(x_1)-\psi(x_2)=I(x_1)-I(x_2)$$
and put
$$\tilde I(x)=I(x_3-x)-\psi(x_3-x)$$
for all $x\in X$. It is easy to check that
$$\tilde I(\lambda(x_1-x_2))=\tilde I((1-\lambda)(x_2-x_1))<\tilde I(0)\ .\eqno{(7)}$$
Fix a
convex set $Y\subseteq X^*$ dense in $X^*$ and put
$$\tilde Y=-Y-\psi\ .$$
Hence, $\tilde Y$ is convex and dense in $X^*$ too.  Now, set
$$A=\{\lambda(x_1-x_2),(1-\lambda)(x_2-x_1)\}\ .$$
Clearly, we have
$$X^*\subset E_{\tilde I, A}\ .\eqno{(8)}$$ 
Since $\tilde Y$ is dense in $X^*$, we have
$$\sup_{\varphi\in \tilde Y}\varphi(x)=+\infty$$
for all $x\in X\setminus \{0\}$. Hence, in view of $(7)$ and $(8)$, we can apply
Theorem 4 with  $B=\{0\}$, $I=\tilde I$, $Y=\tilde Y$. Accordingly,
there exists $\tilde\varphi\in Y$ such that the functional $\tilde I-\tilde\varphi-\psi$
has two global minima in $X$, say $u_1\ne u_2$. At this point, it is clear that
$x_3-u_1, x_3-u_2$ are two global minima of the functional $I+\tilde\varphi$,
and the proof is complete.\hfill $\bigtriangleup$
\medskip
REMARK 2. - We remark that Corollary 1 was also obtained very recently in [7] by
means of a radically different proof.\par
\medskip
We now establish the following technical proposition:\par
\medskip
PROPOSITION 1. - {\it If  $V$ is an infinite-dimensional real Banach space
space and if $U\subset V$ is a $\sigma$-compact set, then there exists
a convex cone $C\subset V$ dense in $V$, such that $U\cap C=\emptyset$.}\par
\smallskip
PROOF. This proposition was proved in [3] when $V$ is a Hilbert space ([3], Proposition 2.4).  As in [3], we distinguish two cases.
First, we assume that $V$ is separable. In this case, the proof provided in [3] can be repeated word for word, and so we
omit it. So, assume that $V$ is not separable. Let $\{x_{\gamma}\}_{\gamma\in \Gamma}$ be a Hamel basis of $V$. Set
$$\Lambda=\{\gamma\in \Gamma : x_{\gamma}\not\in \hbox {\rm span}(U)\}$$
and
$$L=\hbox {\rm span}(\{x_{\gamma} : \gamma\in \Lambda\})\ .$$
Clearly, span$(U)$ is separable since $U$ is so. Hence, $\Lambda$ is infinite.
Introduce in $\Lambda$ a total order $\leq$ with no greatest element. Next, for each $\gamma\in \Lambda$, let
$\psi_{\gamma}:L\to {\bf R}$ be a linear functional such that 
$$\psi_{\gamma}(x_{\alpha})=\cases {1 & if $\gamma=\alpha$\cr & \cr 0 & if $\gamma\neq\alpha$\ .\cr}$$
Now, set
$$D=\{x\in L : \exists \beta\in \Lambda : \psi_{\beta}(x)>0\hskip 5pt 
\hbox {\rm and}\hskip 5pt \psi_{\gamma}(x)=0\hskip 5pt \forall \gamma>\beta\}\ .$$
Of course, $D$ is a convex cone.
Fix $x\in L$. So, there is a finite set $I\subset\Lambda$ such that $x=\sum_{\gamma\in I}\psi_{\gamma}(x)x_{\gamma}$.
Now, fix $\beta\in\Lambda$ so that $\beta>\max I$. For each $n\in {\bf N}$, put
$$y_n=x+{{1}\over {n}}x_{\beta}\ .$$
Clearly, $\psi_{\beta}(y_n)={{1}\over {n}}$ and $\psi_{\gamma}(y_n)=0$ for all $\gamma>\beta$. Hence, $y_n\in D$.
Since $\lim_{n\to \infty}y_n=x$, we infer that $D$ is dense in $L$. At this point, it is immediate to check that the
set $D+\hbox {\rm span}(U)$ is a convex cone, dense in $V$, which does not meet $U$.\hfill
$\bigtriangleup$\par
\medskip
{\it Proof of Theorem 1.} Let us prove that $S_T$ is not $\sigma$-compact.
 Arguing by contradiction, assume the contrary. Then, by Theorem A,
for each $\varphi\in X^*\setminus T(S_T)$, the equation
$$T(x)=\varphi$$
has a unique solution in $X$.
Moreover, since $T$ is continuous, $T(S_T)$ is
 $\sigma$-compact too. Therefore, in view of Proposition 1, there is a convex set $Y\subset X^*$
dense in $X^*$, such that $T(S_T)\cap Y=\emptyset$. On the other hand, since $T$ is not monotone, the functional
$J_T$ is not convex and so, thanks to Corollary 1, there is
$\tilde\varphi\in Y$ such that the functional $J_T-\tilde\varphi$ has at least two global minima in $X$ which
are therefore solutions of the equation
$$T(x)=\tilde\varphi\ ,$$
a contradiction. 
Now, let us prove that $T(S_T)$ is not $\sigma$-compact. Arguing by contradiction, assume the contrary.
 Consequently, since $T$ is proper ([6], Theorem 1),
 $T^{-1}(T(S_{T}))$ would
be $\sigma$-compact.  But then,  since $S_T$ is closed and $S_T\subseteq T^{-1}(T(S_T))$,  $S_T$ would be $\sigma$-compact,
a contradiction. The proof is complete.
 \hfill $\bigtriangleup$     \par
\medskip
{\it Proof of Theorem 2.}
By Theorem 1, the set $S_{T}$ is not $\sigma$-compact.
Now, observe that if
$x\in X\setminus {(\tilde S_{T_{|(X\setminus B)}}\cup B)}$, then,
by the inverse function
theorem, $T$ is a local homeomorphism at $x$, and so $x\not\in
S_{T}$. Hence, we have
$$S_{T}\subseteq \tilde S_{T_{|(X\setminus B)}}\cup B\ .$$
We then infer that $\tilde S_{T_{|(X\setminus B)}}$ is not
$\sigma$-compact since, otherwise,
$\tilde S_{T_{|(X\setminus B)}}\cup B$ would be so, and hence
also $S_{T}$ would be
$\sigma$-compact being closed. Finally, the fact that
$T(\tilde S_{T_{|(X\setminus B)}})$ is not $\sigma$-compact follows
as in the final part of the proof of Theorem 1, taking into account that
$\tilde S_{T_{|(X\setminus B)}}$ is closed in the open set $X\setminus B$ and so it turns out to be the union of an at most countable
family of closed sets. \hfill $\bigtriangleup$\par
\medskip
{\it Proof of Theorem 3.} Clearly,  since $X$ is a Hilbert space, we are identifying $X^*$ to $X$.
 Let us prove that $(i)\to (iii)$. So, assume $(i)$.
Since $J_T$ is not convex, by a classical characterization ([8], Theorem 2.1.11), the set
$\Gamma$ is non-empty. Fix $\mu\in ]\tilde\mu,0[$. For each $x\in X$, put
$$I_{\mu}(x):={{1}\over {2}}\|x\|^2-{{1}\over {\mu}} J_T(x)\ .$$
Clearly, for some $(x,y)\in \Gamma$, we have
$$\left \langle y-{{1}\over {\mu}}T'(x)(y), y\right \rangle <0$$
and so, since 
$$I''_{\mu}(x)(y)= y-{{1}\over {\mu}}T'(x)(y)\ ,$$
the above recalled characterization implies that the functional $I_{\mu}$ is not convex.
Since $T$ is compact, on the one hand, $J_T$ is sequentially
weakly continuous ([10], Corollary 41.9) and, on the other hand, in view of $(2)$ the operator
$I_{\mu}'$ (recall that $-{{1}\over {\mu}}>\lambda_0$) is proper ([9], Example 4.43).
 The compactness of $T$ also implies that, for each
$x\in X$, the operator $T'(x)$ is compact ([9], Proposition 7.33) and so, for each $\lambda\in {\bf R}$,
 the operator $y\to y+\lambda T'(x)(y)$ is Fredholm of index zero ([9], Example 8.16).
Therefore, the operator
 $I_{\mu}'$ is non-monotone, proper and Fredholm of index zero.
Clearly, by $(1)$, the functional $x\to I_{\mu}(x)+\langle z,x\rangle$ is coercive for all $z\in X$. Then, in view of 
Corollary 1, the operator $I_{\mu}'$ is not injective.  At this point, we can apply Theorem B
to infer that the set $\tilde S_{I_{\mu}'}$ contains an accumulation point. Finally, notice that
$$\tilde S_{I_{\mu}'}=A_{\mu}\ ,$$
and $(iii)$ follows.
The implication $(iii)\to (ii)$ is trivial. Finally, the implication $(ii)\to (i)$ is
provided by Theorem 2.1.11 of [8] again.\hfill $\bigtriangleup$
\medskip
REMARK 3. - Some applications of the above results to weighted  eigenvalue problems (which cannot be obtained by means of the results in [3]) are
presented in [5].
\bigskip
\bigskip
\noindent
{\bf Acknowledgements.} The author has been supported by the Gruppo Nazionale per l'Analisi Matematica, la Probabilit\`a e le loro Applicazioni (GNAMPA) of the Istituto Nazionale di Alta Matematica (INdAM). Also, the author thanks the referees and Professor Reich for their remarks.
 \vfill\eject
\centerline {\bf References}\par
\bigskip
\bigskip
\noindent
[1]\hskip 5pt R. A. PLASTOCK, {\it Nonlinear Fredholm maps of index zero and
their singularities}, Proc. Amer. Math. Soc., {\bf 68} (1978), 317-322.\par
\smallskip
\noindent
[2]\hskip 5pt B. RICCERI, {\it A further improvement of a minimax theorem of
Borenshtein and Shul'man}, J. Nonlinear Convex Anal., {\bf 2} (2001),
279-283.\par
\smallskip
\noindent
[3]\hskip 5pt B. RICCERI, {\it On the singular set of certain potential operators in
Hilbert spaces}, in Progr. Nonlinear Differential Equations Appl., {\bf 75},
377-391, Birkh\"auser, 2007. \par
\smallskip
\noindent
[4]\hskip 5pt B. RICCERI, {\it A strict minimax inequality criterion and some of its consequences}, Positivity, 
{\bf 16} (2012), 455-470.\par
\smallskip
\noindent
[5]\hskip 5pt B. RICCERI, {\it Weights sharing the same eigenvalue}, 
arXiv:1407.3439v1 [math.AP], submitted.\par
\smallskip
\noindent
[6]\hskip 5pt R. S. SADYRKHANOV, {\it On infinite                
dimensional features of proper and closed mappings}, Proc. Amer.
Math. Soc., {\bf 98} (1986), 643--658.\par
\smallskip
\noindent
[7]\hskip 5pt J. SAINT RAYMOND, {\it Characterizing convex functions by variational properties},
 J. Nonlinear Convex Anal., {\bf 14} (2013), 253-262. \par
\smallskip
\noindent
[8]\hskip 5pt C. Z\u{A}LINESCU, {\it Convex analysis in general vector spaces}, World Scientific, 2002.\par
\smallskip
\noindent
[9]\hskip 5pt E. ZEIDLER, {\it Nonlinear functional analysis and its
applications}, vol. I, Springer-Verlag, 1986.\par
\smallskip
\noindent
[10]\hskip 5pt E. ZEIDLER, {\it Nonlinear functional analysis and its
applications}, vol. III, Springer-Verlag, 1985.\par
\bigskip
\bigskip
\bigskip
\bigskip
Department of Mathematics\par
University of Catania\par
Viale A. Doria 6\par
95125 Catania\par
Italy\par
{\it e-mail address}: ricceri@dmi.unict.it

\bye